\documentclass[journal]{IEEEtran}
\usepackage{color}
\usepackage[pdftex]{graphicx}
\DeclareGraphicsExtensions{.pdf,.jpeg,.png}
\graphicspath{{./figures/}}

\usepackage[cmex10]{amsmath}
\usepackage{amssymb, latexsym}
\usepackage{amsfonts}
\usepackage{amsthm}
\usepackage{mathrsfs}    

\newcommand{\be}{\begin{equation}}
\newcommand{\ee}{\end{equation}}
\newcommand{\bea}{\begin{eqnarray}}
\newcommand{\eea}{\end{eqnarray}}
\def\bal#1\eal{\begin{align}#1\end{align}}
\renewcommand{\Vec}[1]{\boldsymbol{#1}}
\newcommand{\Mat}[1]{\boldsymbol{#1}}
\newcommand{\Exp}[1]{\mathrm{e}^{#1}}
\newcommand{\trace}[1]{\mathop{\mathrm{tr}}\left({#1}\right)}

\newcounter{saveeqn}

\newcommand{\url}[1]{{\tt\footnotesize {#1}}}

\newcommand{\Vecgamma}{\Vec{\gamma}\hspace*{-1.3ex}\Vec{\gamma}}

\begin{document}

\title{Multi Snapshot Sparse Bayesian Learning for DOA Estimation}

\author{
Peter Gerstoft \IEEEmembership{Member, IEEE}
and
Christoph F. Mecklenbr\"auker \IEEEmembership{Senior Member,~IEEE,}
and 
Angeliki Xenaki 

\thanks{Supported by 
the Office of Naval Research, Grant Nos. N00014-1110439 and N00014-1310632 (MURI), and FTW Austria's ``Compressed channel state information feedback for time-variant MIMO channels''.}
\thanks{P. Gerstoft is with University of California San Diego, La Jolla, CA.
}
\thanks{C. F. Mecklenbr\"auker is with TU Wien, 
1040 Vienna, Austria}
\thanks{A. Xenaki is with Acoustic Technology, Technical University of Denmark.
}}
\maketitle

\begin{abstract}
\textcolor{red}{\today} The directions of arrival (DOA) of plane waves are estimated from multi-snapshot sensor array data using Sparse Bayesian Learning (SBL).
The prior  source amplitudes  is  assumed 
 independent zero-mean complex Gaussian distributed with hyperparameters the unknown variances (i.e. the source powers). For a complex Gaussian likelihood with hyperparameter the unknown noise variance, the corresponding Gaussian posterior distribution is derived.
For a given number of DOAs, the hyperparameters are automatically selected by maximizing the evidence and promote sparse DOA estimates. The  SBL scheme for DOA estimation is discussed  and evaluated competitively against LASSO ($\ell_1$-regularization), conventional beamforming, and MUSIC.
\end{abstract}
%
\IEEEpeerreviewmaketitle

\section{Introduction}

In direction of arrival (DOA) estimation, compressive beamforming, i.e.\ sparse processing, achieves high-resolution acoustic imaging and reliable DOA estimation even with a single snapshot\cite{MalioutovDOA:2005,EdelmannCSDOA:2011,Mecklenbrauker:2013,FortunatiDOA:2014,XenakiCS:2014,Xenaki:2015}, outperforming traditional methods\cite{VanTreesBook}. 

Multiple measurement vector (MMV, or multiple snapshots) compressive beamforming offers several benefits over established high-resolution DOA estimators which utilize the data covariance matrix\cite{MalioutovDOA:2005,XenakiCS:2014,Wipf2007, Gerstoft2015}: 
1) It handles partially coherent arrivals.
2) It can be formulated with any number of snapshots in contrast to eigenvalue based beamformers. 
3) Its flexibility in formulation enables extensions to sequential processing, and online algorithms \cite{Mecklenbrauker:2013}. 
4) It achieves higher resolution than MUSIC, even in scenarios that favor these classical high-resolution methods \cite{Gerstoft2015}.

We solve the MMV problem in the sparse Bayesian learning (SBL) framework\cite{Wipf2007} and use the maximum-a-posteriori (MAP) estimate for DOA reconstruction. We assume complex Gaussian distributions with unknown variances (hyperparameters) both for the likelihood and as prior information for the  source amplitudes. Hence, the corresponding posterior distribution is also Gaussian. 
%
To determine the hyperparameters, we maximize a Type-II  likelihood (evidence)  for Gaussian signals hidden in Gaussian noise.
This has been solved with a Minimization-majorization based technique\cite{Stoica2012} and with expectation maximization (EM) \cite{Wipf2007,Wipf2007beam,Wipf2004,Zhang2011,Liu2012, Zhang2016}. Instead, we estimate the hyperparameters directly from the likelihood derivatives using stochastic maximum likelihood\cite{Boehme1985,Jaffer1988,Stoica-Nehorai1995}.

We propose a SBL algorithm for MMV DOA estimation which, given the number of sources, automatically estimates the set of DOAs corresponding to non-zero source power from all potential DOAs.
This provides a sparse signal estimate similar to LASSO\cite{TibshiraniLasso1996,Gerstoft2015}. Posing the  problem this way, the  estimated number of parameters is  independent of snapshots, while the accuracy  improves with the number of snapshots. 

\section{Array data model and problem formulation}

Let $\Mat{X}\!\!=\!\![\Vec{x}_1,\ldots,\Vec{x}_L]\!\!\in\!\!\mathbb{C}^{M\times L}$ be the complex source amplitudes, $x_{ml}$ with $m\in[1,\cdots, M]$ and $l\in[1,\cdots, L]$, at $M$ DOAs (e.g. $\theta_m=-90^{\circ}+\frac{m-1}{M}180^{\circ}$) and L snapshots at frequency $\omega$. We observe narrowband waves on  $N$ sensors for $L$ snapshots
$\Mat{Y}=[\Vec{y}_1,\ldots,\Vec{y}_L]\in\mathbb{C}^{N\times L}$.
A linear regression model relates the  array data 
$ \Mat{Y} $  to the source amplitudes  $\Mat{X} $, 
\be
   \Mat{Y} = \Mat{A}\Mat{X} + \Mat{N}~.  \label{eq:linear-model3}
 \ee
The transfer matrix $ \Mat{A}=[\Vec{a}_1 \ldots,\Vec{a}_M]\in\mathbb{C}^{N\times M}$ contains the array
steering vectors for all hypothetical DOAs as columns, with the $nm$th element 
$\mathrm{e}^{-\mathrm{j}(n-1)\frac{\omega d}{c}\sin\theta_m}$ ($d$ is the element spacing and $c$ the sound speed). 
The  additive noise  $\Vec{N}\!\!\in\!\!\mathbb{C}^{N\times L}$ is assumed independent across sensors and snapshots, with each element following a complex Gaussian ${\cal CN}(0, \sigma^2)$.

We assume $M\gg N$ and thus (\ref{eq:linear-model3}) is underdetermined. In the presence of few stationary sources, the source vector $\Vec{x}_l$ is  $K$-sparse with $K\ll M$.
We define the $l$th active set
\be
     \mathcal{M}_l = \{m\in\mathbb{N}| x_{ml}\ne0\} = \{m_1,\,m_2,\ldots,\,m_K \}~,
\ee
and assume $  \mathcal{M}_l =   \mathcal{M}$  is constant across snapshots $l$.
Also,  we define $\Mat{A}_{\mathcal{M}}\in\mathbb{C}^{N\times K}$  which contains only the $K$
``active'' columns of $\Mat{A}$.
The $\|\cdot\|_p$  denotes the vector $p$-norm and $\|\cdot\|_{\mathcal{F}}$ the matrix Frobenius norm.


\section{Bayesian formulation}
\label{sec:Bayes}

Using Bayesian inference to solve the linear problem \eqref{eq:linear-model3} involves determining the posterior distribution of the complex source amplitudes $\Mat{X}$ from the  likelihood and a prior model.

\subsection{Likelihood}
Assuming the additive noise  \eqref{eq:linear-model3} complex Gaussian the data likelihood, i.e., the conditional probability density function (pdf)
for the single-frequency observations $\Vec{Y}$ given the sources $\Vec{X}$, is  complex Gaussian with noise variance $\sigma^2$.
\begin{align}\label{eq:likelihood}
     p(\Vec{Y} | \Vec{X}; \sigma^2) 
   &=  \frac{ \exp\left( - \frac{1}{\sigma^2} \| \Vec{Y} - \Mat{A}\Mat{X}\|_{\mathcal{F}}^2 \right)}{(\pi\sigma^2)^{NL}} .
\end{align}

\subsection{Prior}

We assume that the complex source amplitudes $x_{ml}$ are independent both across snapshots and across DOAs and follow a zero-mean complex Gaussian distribution
with DOA-dependent variance $\gamma_m \in \Vecgamma=[\gamma_1,\ldots,\gamma_M]^T$,
\begin{align}
 \label{eq:xprior}
    p_m(x_{ml}; \gamma_m) &= \left\{ \begin{array}{ll}
        \delta(x_{ml}), & \text{for } \gamma_m = 0, \\
        \frac{1}{\pi\gamma_m} \mathrm{e}^{-|x_{ml}|^2/\gamma_m}, & \text{for } \gamma_m > 0
        \end{array} \right.
\\
\label{eq:tipping-like}
p(\Mat{X};\Vecgamma) 
&=  \prod_{l=1}^L \prod_{m=1}^M p_m(x_{ml}; \gamma_m)=\prod_{l=1}^L \mathcal{CN}(\Vec{0},\Mat{\Gamma}),
\end{align}
i.e., the source vector $\Vec{x}_{l}$ at each snapshot $l\in[1,\cdots, L]$ has a multivariate Gaussian distribution with potentially singular covariance matrix, 
\be
\Mat{\Gamma} = \mathop{\mathrm{diag}}(\Vecgamma)= \mathop{\mathsf{E}}\left[   \Vec{x}_l\Vec{x}_l^H; \Vecgamma\right], 
\ee
as $\mathop{\mathrm{rank}}(\Mat{\Gamma}) = \mathop{\mathrm{card}}(\mathcal{M}) = K \le M$. Note that the diagonal elements of $\Mat{\Gamma} $, i.e., the hyperparameters $\Vecgamma\ge\bf{0}$,  represent source powers. When the variance $\gamma_{m}=0$, then $x_{ml}=0$ with probability 1. 
The sparsity of the model is thus controlled with the hyperparameters $\Vecgamma$.

\subsection{Posterior}
Given the likelihood for the array observations $\Mat{Y}$ \eqref{eq:likelihood} and the prior \eqref{eq:tipping-like}, the posterior pdf for the source amplitudes $\Mat{X}$ can be found using Bayes rule conditioned on 
$\Vecgamma,\sigma^2$,
\be\label{eq:Tipping(31)}
    p(\Mat{X} | \Mat{Y}; \Vecgamma, \sigma^2) 
    \equiv
    \frac{p(\Mat{Y}|\Mat{X}; \sigma^2 )p(\Mat{X}; \Vecgamma)}
            {p(\Mat{Y} ; \Vecgamma,\sigma^2)}.
\ee
The denominator $p(\Mat{Y} ; \Vecgamma,\sigma^2)$ is the evidence term, i.e., the marginal distribution for the data, which for a given $\Vecgamma,\sigma^2$ is a normalization factor  and is neglected at first,
\begin{align}
   p(\Mat{X} | &\Mat{Y}; \Vecgamma,\sigma^2) 
   \propto   p(\Mat{Y}|\Mat{X};\sigma^2) p(\Mat{X}; \Vecgamma) \label{eq:product} \\
       &\propto \frac{\Exp{-\trace{(\Mat{X}-\Mat{\mu}_{\Mat{X}})^H\Mat{\Sigma}_x^{-1}(\Mat{X}-\Mat{\mu}_{\Mat{X}})}}}{(\pi^{N}\det\Mat{\Sigma}_x)^L} 
        ={\cal CN}(\Mat{\mu}_{\Mat{X}}, \Vec{\Sigma_{\Vec{x}}}) ~.       
\end{align}
As both $p(\Mat{Y}|\Mat{X};\sigma^2)$  in \eqref{eq:likelihood} and $p(\Mat{X}; \Vecgamma)$ in \eqref{eq:tipping-like} are Gaussians, their product (\ref{eq:product}) is Gaussian with posterior mean $\Mat{\mu}_{\Mat{X}}$ and covariance $\Vec{\Sigma_{\Vec{x}}}$,
\begin{align}
    \Mat{\mu}_{\Mat{X}}   &= \mathop{\mathsf{E}}\{ \Mat{X} | \Mat{Y}; \Vecgamma,\sigma^2  \} = 
    \Mat\Gamma \Mat{A}^H\Sigma_{\Vec{y}}^{-1}\Mat{Y} ,\label{eq:Wipf-and-Rao-2007} \\ 
    \Mat{\Sigma}_{\Vec{x}} &=  \mathop{\mathsf{E}}\{ (\Vec{x}_l -   \Mat{\mu}_{\Vec{x}_l}  )(\Vec{x}_l -   \Mat{\mu}_{\Vec{x}_l}  )^H | \Vec{Y}; \Vecgamma,\sigma^2  \} \nonumber \\
    &=
    \left( \frac{1}{\sigma^2} \Mat{A}^H \Mat{A} + \Mat{\Gamma}^{-1} \right)^{-1}  
   =\Mat{\Gamma} -\Mat{\Gamma} \Mat{A}^H \Sigma_{\Vec{y}}^{-1} \Mat{A} \Mat{\Gamma} , \label{eq:Sigma_x}
\end{align}
where the array data covariance  $\Mat\Sigma_{\Vec{y}}$ and its inverse are derived from \eqref{eq:linear-model3} and using the matrix inversion lemma
\begin{align}
  \Mat\Sigma_{\Vec{y}}&=  \mathop{\mathsf{E}}\{ \Vec{y}_l \Vec{y}_l^H \}  = \sigma^{2}\Mat{I}_N +  \Mat{A} \Mat{\Gamma}\Mat{A}^H  , \label{eq:array-covariance} \\
 \Mat\Sigma_{\Vec{y}}^{-1}&= \sigma^{-2}\Mat{I}_N-\sigma^{-2}\Mat{A}
   \left( \frac{1}{\sigma^2} \Mat{A}^H \Mat{A} + \Mat{\Gamma}^{-1} \right)^{-1} 
   \Mat{A}^H\sigma^{-2} \nonumber \\
 &= \sigma^{-2}\Mat{I}_N-\sigma^{-2}\Mat{A}\Mat{\Sigma}_{\Vec{x}}\Mat{A}^H\sigma^{-2}.
    \label{eq:Tipping(33)} 
\end{align}
If $\Vecgamma$ and $\sigma^2$ are known then the MAP estimate is the posterior mean,
\be
\hat{\Mat X}^{\mathrm{MAP}}=\Vec{\mu}_{\Mat X}=\Mat\Gamma \Mat{A}^H\Mat{\Sigma}_{\Vec{y}}^{-1}\Mat{Y} ~.
\label{eq:x-map-estimate}
\ee	
The diagonal elements of $\Mat\Gamma $ control the row-sparsity of $\hat{\Mat X}^{\mathrm{MAP}}$ as for $\gamma_m=0$ 
the corresponding $m$th row of $ \hat{\Mat X}^{\mathrm{MAP}}$ becomes $\Vec{0}^T$.
Thus, the active set $\mathcal{M}$ is equivalently defined  by 
\be
     \mathcal{M} = \{m\in\mathbb{N} | \gamma_m > 0 \} \label{eq:active-set-gamma}~.
\ee

\subsection{Evidence}
The hyperparameters $\Vecgamma,\sigma^2$  in (\ref{eq:Wipf-and-Rao-2007}--\ref{eq:Tipping(33)}) are estimated by a type-II maximum likelihood, i.e., by maximizing the evidence which was treated as constant in \eqref{eq:product}. The evidence is the product of the likelihood (\ref{eq:likelihood}) and the prior (\ref{eq:tipping-like}) integrated over the  complex  source amplitudes $\Mat{X}$, 
\begin{align}
p(\Mat{Y} ; \Vecgamma, \sigma^2) &= \int_{{\mathbb R}^{2ML}}
    p(\Mat{Y} | \Mat{X}; \sigma^2)
    p(\Mat{X} ; \Vecgamma)  \, {\rm d}\Mat{X}
                \nonumber \\
    &= \frac{\Exp{-\trace{\Mat{Y}^H \Mat\Sigma_{\Vec{y}}^{-1}\ \Mat{Y}}}}{(\pi^{N}\det\Mat\Sigma_{\Vec{y}})^L}~,
\end{align}
%
where  ${\rm d}\Mat{X}=\prod_{l=1}^L\prod_{m=1}^M 
 \mathop{\mathrm{Re}}(\mathrm{d}X_{ml}) \mathop{\mathrm{Im}}(\mathrm{d}X_{ml})$, and
$\Mat\Sigma_{\Vec{y}}$ is the data covariance \eqref{eq:array-covariance}. 
The $L$-snapshot marginal log-likelihood becomes
\begin{align}
\log p(\Mat{Y} ; \Vecgamma, \sigma^2)
    &\propto -\trace{\Mat{Y}^H \Mat\Sigma_{\Vec{y}}^{-1} \Mat{Y}} - L \log\det\Mat\Sigma_{\Vec{y}} \nonumber \\
    &\propto -\trace{\Mat\Sigma_{\Vec{y}}^{-1} \Mat{S}_{\Vec{y}}} - \log\det\Mat\Sigma_{\Vec{y}},
    \label{eq:evidence-L}
\end{align}
where we define the data sample covariance matrix,
\be
     \Mat{S}_{\Vec{y}} =  \Mat{YY}^H/L.
\ee
Note that \eqref{eq:evidence-L} does not involve the inverse of $\Mat{S}_{\Vec{y}}$ hence it works well even for few snapshots (small L). 

The  hyperparameter estimates $\hat{\Vecgamma},\hat{\sigma}^2$ are obtained by maximizing the evidence,
\be
\label{eq:TypeIIMaxLikelihood}
   (\hat{\Vecgamma},\;\hat{\sigma}^2) = \mathop{\arg\max}_{\Vecgamma\ge0,\;\sigma^2>0} \log  p(\Mat{Y} ; \Vecgamma, {\sigma^2}).
\ee
The maximization is carried out iteratively using derivatives of the evidence for $\Vecgamma$ (see Sec. \ref{sec:src-pwr-estim}) as well as conventional noise estimates (see Sec. \ref{se:noise}) as explained in Sec. \ref{eq:sbl-algo}.

\subsection{Source power estimation (hyperparameters $\Vecgamma$)}
\label{sec:src-pwr-estim}
We impose the diagonal structure $\Mat{{\Gamma}}=\mathop{\mathrm{diag}}
(\Vec{{\gamma}})$, in agreement with \eqref{eq:tipping-like}, and form derivatives of \eqref{eq:evidence-L} with respect to the diagonal elements $\gamma_m$,
cf. \cite{Boehme1986}.
Using
\begin{align}
&\frac{\partial \Mat{\Sigma}_{\Vec{y}}^{-1}}{\partial \gamma_m} = -\Mat{\Sigma}_{\Vec{y}}^{-1} \, \frac{\partial \Mat{\Sigma}_{\Vec{y}}}{\partial \gamma_m}  \, \Mat{\Sigma}_{\Vec{y}}^{-1}
  = -\Mat{\Sigma}_{\Vec{y}}^{-1} \Vec{a}_m^{}\Vec{a}_m^H \Mat{\Sigma}_{\Vec{y}}^{-1}, \label{eq:deriv-inverse}
\\
& \frac{\partial\log\det( \Mat{\Sigma}_{\Vec{y}})}{\partial\gamma_m} 
 = \trace{\Mat{\Sigma}_{\Vec{y}}^{-1} \frac{\partial \Mat{\Sigma}_{\Vec{y}}}{\partial\gamma_m}} 
 =   \Vec{a}_m^H   \Mat{\Sigma}_{\Vec{y}}^{-1} \Vec{a}_m 
 \label{eq:deriv-log-det},
\end{align}
the derivative of \eqref{eq:evidence-L} is
\begin{align}
 \frac{\partial\log p(\Mat{Y};\Vecgamma,\sigma^2)}{\partial\gamma_m} 
    &= \frac{1}{\gamma_m^{2}L} \| \Vec{\mu}_m  \|_2^2
    -  \Vec{a}_m^H \Mat{\Sigma}_{\Vec{y}}^{-1} \Vec{a}_m,
    \label{eq:score-gamma}
\end{align}
where $\Vec{\mu}_m=\gamma_m\Vec{a}_m^H\Mat{\Sigma}_{\Vec{y}}^{-1}\Mat{Y}$ is the $m$th row of $\Mat{\mu}_{\Mat{X}}$ in \eqref{eq:Wipf-and-Rao-2007}. 
Assuming $\Vec{\mu}_m$ given (from previous iterations or initialization) and forcing \eqref{eq:score-gamma} to zero gives the $\gamma_m$ update (\ref{eq:gamma-update}):
\be
    \gamma_m^{\rm new} = {\frac{1}{\sqrt{L}}\|\Vec{\mu}_m\|_2}/{\sqrt{\Vec{a}_m^H\Mat{\Sigma}_{\Vec{y}}^{-1}\Vec{a}_m}}.
  \tag{SBL1}
   \label{eq:gamma-update}
\ee

When the sample data covariance 
$\Mat{S}_{\Vec{y}}$ is positive definite (i.e.\ usally when $L\ge 2N$) we can replace  $\Mat{\Sigma}_{\Vec{y}}^{-1}$ in \eqref{eq:gamma-update} 
with  $ \Mat{S}_{\Vec{y}}^{-1}$  [see \eqref{eq:necessary-condition}]
\be
    \gamma_m^{\rm new} = 
  {\frac{1}{\sqrt{L}}\|\Vec{\mu}_m\|_2}/{\sqrt{\Vec{a}_m^H\Mat{S}_y^{-1}\Vec{a}_m}}.
\tag{SBL}
\label{eq:improved-update}
\ee
The \ref{eq:improved-update} estimate tends to converge faster as the denominator does not change during iterations.

Wipf and Rao (\cite{Wipf2007}: Eq.(18)) followed the EM approach to estimate the update \ref{eq:gamma-update-EM}:
\be
\gamma_m^{\rm new}={\frac{1}{L}\|\Vec{\mu}_m\|_2^2}+(\Mat{\Sigma}_{\Vec{x}})_{mm} .
   \tag{M-SBL}
   \label{eq:gamma-update-EM}
\ee

The sequence of parameter estimates in the EM iteration has been proven to converge \cite{dempster1977}. However, the convergence is only guaranteed towards a \emph{local} optimum of the marginal log-likelihood \eqref{eq:evidence-L}. As shown in Sec. \ref{se:example} all the update rules \eqref{eq:gamma-update}--\eqref{eq:gamma-update-EM} converge provided 
$| {\partial \gamma_m^{\rm new}}/{\partial \gamma_{m}}|<1$. 

 \begin{table}
\begin{tabular}{rlc} \hline\hline
0& Given: $\Mat{A}\in\mathbb{C}^{N\times M}$,\; $\Mat{Y}\in\mathbb{C}^{N \times L}$,\; $K=3$ & \\ 
    & Initialize, here: $\sigma_0^2=0.1, \gamma_0=1, \epsilon_{\min}=0.001, j_{\max}=500$ \hspace*{-30ex}\\ \hline
1& initialize $j=0$, $\sigma^2=\sigma_0^2$, $\Vecgamma = \gamma_0$ & \\
2& while $(\epsilon > \epsilon_{\min})$ and $(j<j_{\max})$  & \\
3&  \hspace*{2ex} $j=j+1, \quad$   $\Vecgamma^{\mathrm{old}} = \Vecgamma^{\mathrm{new}} , \quad \Mat{\Gamma} = \mathop{\mathrm{diag}}(\Vecgamma^{\mathrm{new}} )$ \hspace*{-30ex} & \\
4&  \hspace*{2ex}  $\Mat{\Sigma}_{\Vec{y}} = \sigma^2\Mat{I}_N + \Mat{A} \Mat{\Gamma} \Mat{A}^H$ & \eqref{eq:array-covariance} \\ 
5&  \hspace*{2ex} $\Vec{\mu}_m =\gamma_m\Vec{a}_m^H\Mat{\Sigma}_{\Vec{y}}^{-1}\Mat{Y}$ & \eqref{eq:Wipf-and-Rao-2007} \\   
6&  \hspace*{2ex}  $\gamma_m^{\mathrm{new}} = \left\{ \begin{array}{l}\frac{1}{\sqrt{L}} \|\Vec{\mu}_m\|_2 \Big/ \sqrt{ \Vec{a}_m^H \Mat{S}_{\Vec{y}}^{-1}\Vec{a}_m} \\
\frac{1}{\sqrt{L}} \|\Vec{\mu}_m\|_2 \Big/ \sqrt{ \Vec{a}_m^H \Mat{\Sigma}_{\Vec{y}}^{-1}\Vec{a}_m} \\
 {\frac{1}{L}\|\Vec{\mu}_m\|_2^2}+(\Mat{\Sigma}_{\Vec{x}})_{mm} 
 \end{array}\right. $  & $\begin{array}{c} \text{(\ref{eq:improved-update})} \\*[0.7ex]
 \text{(\ref{eq:gamma-update})} \\*[0.7ex]
 \text{(\ref{eq:gamma-update-EM})} \end{array}$ \\
7&  \hspace*{2ex}   $ \mathcal{M} \!\!\!= \!\!\!\{m \in \mathbb{N} |\, \mbox{K largest peaks in} \,\Vecgamma \}\!= \!\{ m_1\ldots m_K \}$  \hspace*{-30ex}& \eqref{eq:active-set-gamma}\\
8& \hspace*{2ex}  $\Mat{A}_\mathcal{M}= (a_{m_1},\ldots,a_{m_K})$ & \\
9& \hspace*{2ex}   $({\sigma}^2)^{\rm new} =\frac{1}{N-K}\trace{(\Mat{I}_N-\Mat{\Mat{A}_\mathcal{M}\Mat{A}_\mathcal{M}}^+) \Mat{S}_{\Vec{y}}}$ & \eqref{eq:bar-sigma2} \\
10& \hspace*{2ex}  $\epsilon = \| \Vecgamma^{\mathrm{new}}  -\Vecgamma^{\mathrm{old}} \|_1 / \| \Vecgamma^{\mathrm{old}} \|_1$ &  \eqref{eq:improvement} \\
11& end \\ \hline
12& Output: $\mathcal{M}$, $\Vecgamma^{\mathrm{new}}$, $({\sigma}^2)^{\rm new} $ & \\ \hline\hline
\end{tabular}
\caption{SBL Algorithm: In line 6 choose  \ref{eq:improved-update}, \ref{eq:gamma-update}  or \ref{eq:gamma-update-EM}.}
\label{tab:SBLML2}
\end{table}

\subsection{Noise variance estimation (hyperparameter $\sigma^2$)}
\label{se:noise}
Obtaining a good noise variance estimate is important for fast convergence of the SBL method, as it controls the sharpness of the peaks. For a given set of active DOAs $\cal M$, stochastic maximum likelihood \cite{Liu2012,Boehme1985} provides an asymptotically efficient estimate of  $\sigma^2$.

Let  $\Mat{{\Gamma}}_{\cal M}=\mathop{\mathrm{diag}}(\Vecgamma_{\cal M}^{\rm new})$ be the  covariance matrix of the $K$ active sources obtained above with corresponding active steering matrix $\Mat{A}_{\cal M}$ which maximizes (\ref{eq:evidence-L}). 
The corresponding data covariance matrix is
\be
\Mat{{\Sigma}}_{\Vec{y}} = {\sigma}^2\Mat{I}_N + \Mat{A}_{\cal M}\Mat{{\Gamma}}_{\cal M}\Mat{A}_{\cal M}^H ,
\label{eq:stationary-model-covariance}
\ee
where $\Mat{I}_N$ is the identity matrix of order $N$.
The data covariance models \eqref{eq:array-covariance} and \eqref{eq:stationary-model-covariance} are identical. At the optimal solution $(\Mat{{\Gamma}}_{\cal M}, {\sigma}^2)$,  Jaffer's  
necessary condition (\cite{Jaffer1988}:Eq.(6)) must be satisfied
\be
     \Mat{A}_{\cal M}^H \left(  \Mat{S}_{\Vec{y}} - \Mat{{\Sigma}}_{\Vec{y}} \right) \Mat{A}_{\cal M} = \Mat{0}.
\label{eq:necessary-condition}
\ee
Substituting (\ref{eq:stationary-model-covariance}) into (\ref{eq:necessary-condition}) gives
\be
\label{eq:JafferModified}
     \Mat{A}_{\cal M}^H \left(   \Mat{S}_{\Vec{y}} - {\sigma}^2\Mat{I}_N \right) \Mat{A}_{\cal M} =
      \Mat{A}_{\cal M}^H \Mat{A}_{\cal M}\Mat{{\Gamma}}_{\cal M}\Mat{A}_{\cal M}^H \Mat{A}_{\cal M}.
\ee
Multiplying \eqref{eq:JafferModified} from right and left with the pseudo inverse $\Mat{A}_{\cal M}^+ = (\Mat{A}_{\cal M}^H\Mat{A}_{\cal M})^{-1}\Mat{A}_{\cal M}^H$ and  $\Mat{A}_{\cal M}^{+H}$ respectively and subtracting $\Mat{S}_{\Vec{y}}$ from both sides yields \cite{Boehme1985}
\begin{align}
 ({\sigma}^2)^{\rm new} &= \frac{1}{N-K}\trace{(\Mat{I}_N-\Mat{A}_{\cal M}\Mat{A}_{\cal M}^+)  \Mat{S}_{\Vec{y}}}. \label{eq:bar-sigma2}
\end{align}
This estimate requires $K < N$ and will underestimate the noise for small $L$. 

Several estimators for the noise $\sigma^2$ are proposed based on EM
\cite{Wipf2007,Wipf2004,Zhang2011,Tipping2001,Zhang2014}.
Empirically, neither of these  converge well in our application. 
For a comparative illustration in Sec.\ \ref{se:example} we use the iterative noise  $\sigma^2$ EM estimate in \cite{Zhang2014},
\begin{align}
(\sigma^{2})^{\rm new}\!\!\!\! &= \!\!\!\frac{\frac{1}{L}\|(\Mat{Y}\!\!-\!\!\Mat{A} \Mat{\mu}_{\Mat{X}}) \|^2_{\mathcal{F}} +(\sigma^{2})^{\rm old}(M \!\!-\!\!\sum^M_{i=1}\!\!\!\!\!\frac{(\Sigma_x)_{ii}}{\gamma_i})} {N}~. \label{eq:bar-sigma2-Zhang}
\end{align}

\begin{figure}[tb] 
   \includegraphics[width=0.4\textwidth]{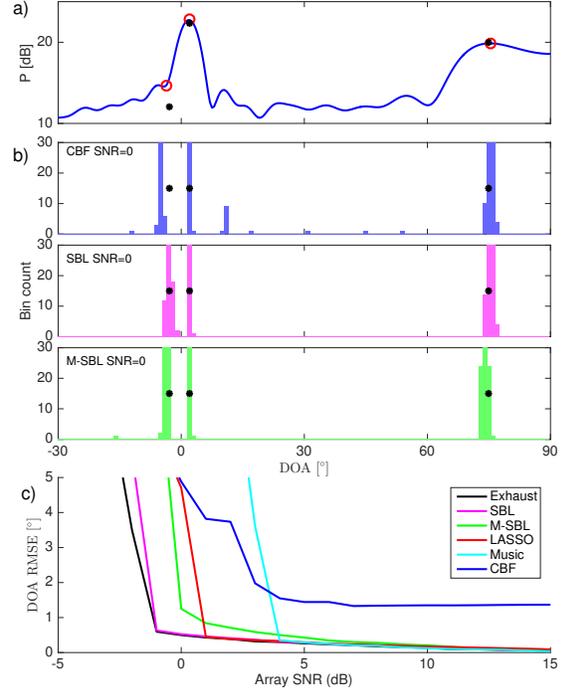}
    \caption{ Multiple $L$=50 snapshot example for  sources at DOAs $[-3,\, 2,\, 75]^\circ$ with magnitudes $[12,\, 22,\, 20]$ dB. a)  spectra for CBF and SBL (o) at SNR=0 dB. b)  CBF, \ref{eq:improved-update}, and \ref{eq:gamma-update-EM} histogram based on 100 Monte Carlo simulations at SNR=0 dB. c) RMSE performance versus array SNR for  exhaustive, \ref{eq:improved-update}, \ref{eq:gamma-update-EM}, LASSO, MUSIC, and CBF. The true source positions ($\bullet$) are indicated in a) and b).    \label{fig:MMV}}
\end{figure}

\subsection{SBL Algorithm}
\label{eq:sbl-algo}
Given the observed $\Mat Y$, we iteratively update ${\Mat \mu}_{\Mat{X}}$  \eqref{eq:Wipf-and-Rao-2007} and ${\Mat \Sigma}_{\Vec{y}}$  \eqref{eq:array-covariance} by using the current $\Vecgamma$.
Either \ref{eq:improved-update}, \ref{eq:gamma-update}, or \ref{eq:gamma-update-EM} can  update  $\gamma_m$ for $m=1,\ldots,M$ and then (\ref{eq:bar-sigma2}) is used to estimate $\sigma^2$. The algorithm is summarized in Table \ref{tab:SBLML2}.

The convergence rate $\epsilon$ measures the relative improvement of the estimated total source power,
\be
\epsilon= \| \Vecgamma^{\rm new}- \Vecgamma^{\rm old} \|_1 \; \Big/ \; 
\| \Vecgamma^{\rm old}\|_1 ~.
\label{eq:improvement}
\ee
The algorithm stops when $\epsilon\le\epsilon_{\min}$ and the output is the
active set $\mathcal{M}$ \eqref{eq:active-set-gamma} from which all relevant source parameter estimates are computed.
%


\begin{figure}[tb] 
   \includegraphics[width=0.45\textwidth]{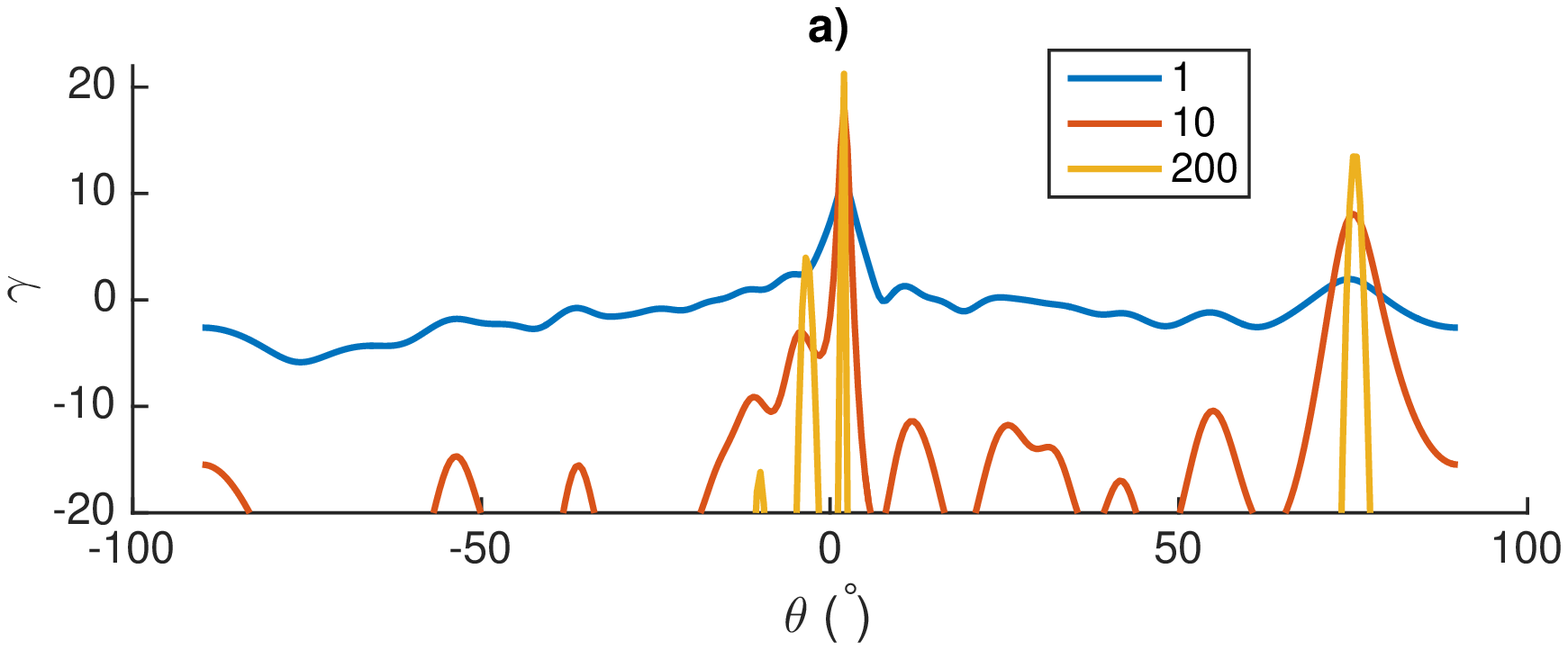}
      \includegraphics[width=0.45\textwidth]{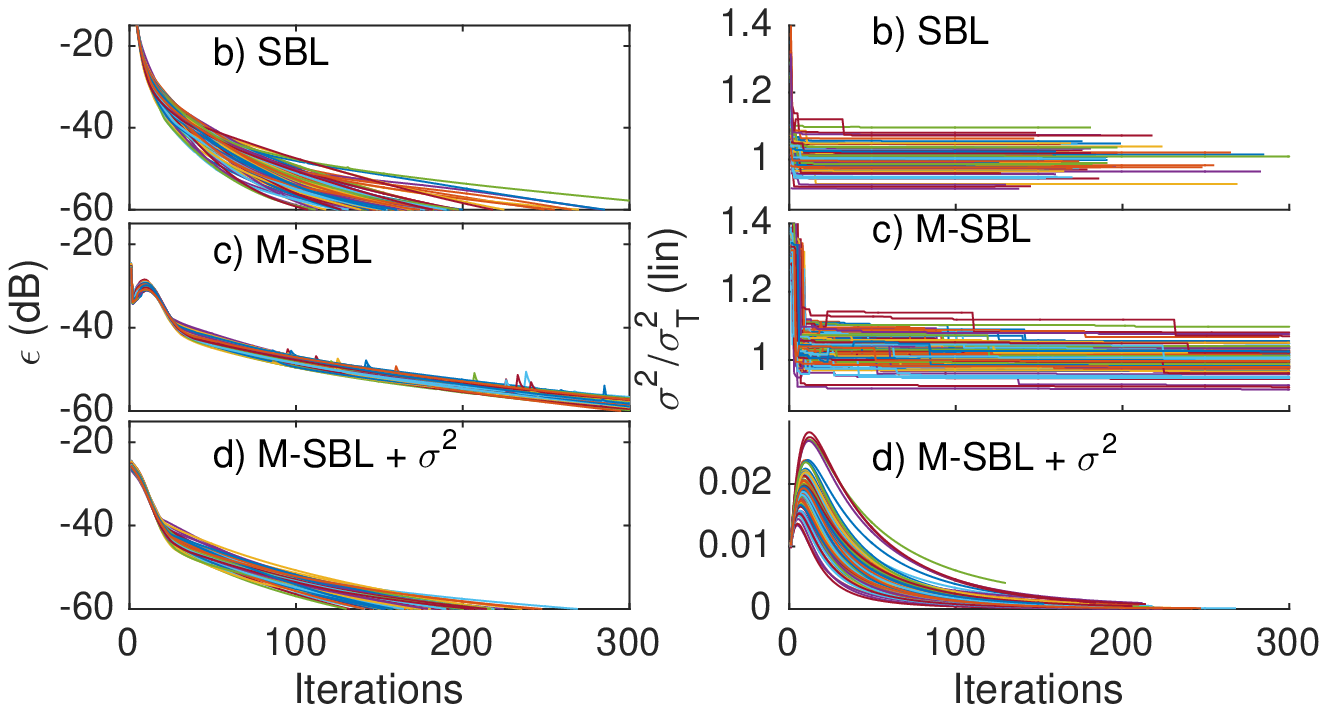}
    \caption{Convergence at SNR=0 dB  with $L$=50. a) $\Vecgamma$ at iteration 1, 10 , 200 for  \ref{eq:improved-update}.
   Convergence of (b)  \ref{eq:improved-update}  and (c  and d) \ref{eq:gamma-update-EM}    for 100 Monte Carlo simulations.
Convergence is shown  for $\epsilon$ (left) and $\sigma^2/\sigma^2_T$ (right).
In (b--c) the noise estimate  $(\sigma^2)^{\rm new}$ is based on \eqref{eq:bar-sigma2} and in d) \eqref{eq:bar-sigma2-Zhang}. } \label{fig:convergence}
\end{figure}

\section{Example}
\label{se:example}

For multiple sources with well separated DOAs and similar magnitudes, conventional beamforming (CBF) and LASSO/SBL methods provide similar DOA estimates. They differ, however, in their behavior whenever two sources are closely spaced. Thus, we examine 3 sources at DOAs $[-3,\, 2,\, 75]^\circ$ with magnitudes $[12,\, 22,\, 20]$ dB\cite{Gerstoft2015}.


We consider an array with N=20 elements and half wavelength intersensor spacing.  The DOAs  are assumed to be on an angular grid [$-90$:$0.5$:$90$]$^{\circ}$, M=361, and  L=50 snapshots are observed. The noise is modeled as iid complex Gaussian, though robustness to array imperfections \cite{Weiss2014} and extreme noise distributions \cite{Ollila2015} can be important. The single-snapshot array signal-to-noise ratio (SNR) is
$
\mathrm{SNR}=10\log_{10}[
{ \mathop{\mathsf{E}}\left\lbrace\lVert\Mat{A}\Vec{x}_l\rVert_{2}^{2}\right\rbrace}/{ \mathop{\mathsf{E}}\left\lbrace\lVert\Vec{n}_l\rVert_{2}^{2}\right\rbrace} ]
\label{eq:ArraySNR}
$. 
Then, for $L$ snapshots the noise power ${\sigma}^2_T$ is
\be
{\sigma}^2_T=\mathop{\mathsf{E}}[ \| \Vec{N} \|_{\mathcal{F}}^2]/L/N = 10^{-{\rm SNR}/10} \; \mathop{\mathsf{E}} \|\Mat{AX}\|^2_{\mathcal{F}}/L/N.
\label{eq:sigma_exp}
\ee 
The estimated $({\sigma}^2)^{\rm new}$  \eqref{eq:bar-sigma2} deviates from  ${\sigma}^2_T$ \eqref{eq:sigma_exp} randomly.


Figure \ref{fig:MMV} compares DOA estimation methods for the simulation. The LASSO solution is found considering multiple snapshots \cite{Gerstoft2015} and programmed in CVX\cite{CVX}. \ref{eq:improved-update} and  \ref{eq:gamma-update-EM} are calculated using the pseudocode on Table \ref{tab:SBLML2}. CBF  suffers from low-resolution and the effect of sidelobes in contrast to sparsity based methods 
as shown in the power spectra in Fig. \ref{fig:MMV}a. 

At array SNR=0 dB  the histogram in Fig. \ref{fig:MMV}b shows that CBF poorly locates the neighboring DOAs at broadside.  \ref{eq:improved-update} and \ref{eq:gamma-update-EM}  localize the sources well. 
The  root mean squared error (RMSE) in Fig. \ref{fig:MMV}c shows that CBF has low resolution  as the main lobe is too broad (see Fig. \ref{fig:MMV}a) and
MUSIC performs well for $\mathrm{SNR}\!\!>\!\!5\,\mathrm{dB}$.
For this case we  include exhaustive search, which defines a lower  performance bound and requires $361!/(3!358!)\!\!=\!\!7.8\!\cdot \!10^6$ evaluations.
 LASSO and the SBL methods perform better than MUSIC and  
 offer similar accuracy to the exhaustive search.
  
We compare the convergence of  \ref{eq:improved-update} and  \ref{eq:gamma-update-EM} at array SNR=0 dB (Fig.\   \ref{fig:convergence}). The spatial spectrum (Fig.\   \ref{fig:convergence}a) shows how the estimate $\Vecgamma$ improves with \ref{eq:improved-update} iterations from initially locating only the main peak to locating also the weaker sources. 
%
\ref{eq:improved-update} exhibits faster convergence  than \ref{eq:gamma-update-EM} to  $\epsilon_{\rm min}=-$60 dB where the algorithm stops (Figs.\ \ref{fig:convergence}b versus \ref{fig:convergence}c). 
\ref{eq:gamma-update-EM} underestimates  $\sigma^2$ significantly when using \eqref{eq:bar-sigma2-Zhang} (Fig.\ \ref{fig:convergence}d).

\begin{figure}[tb] 
 \center\includegraphics[width=0.4\textwidth]{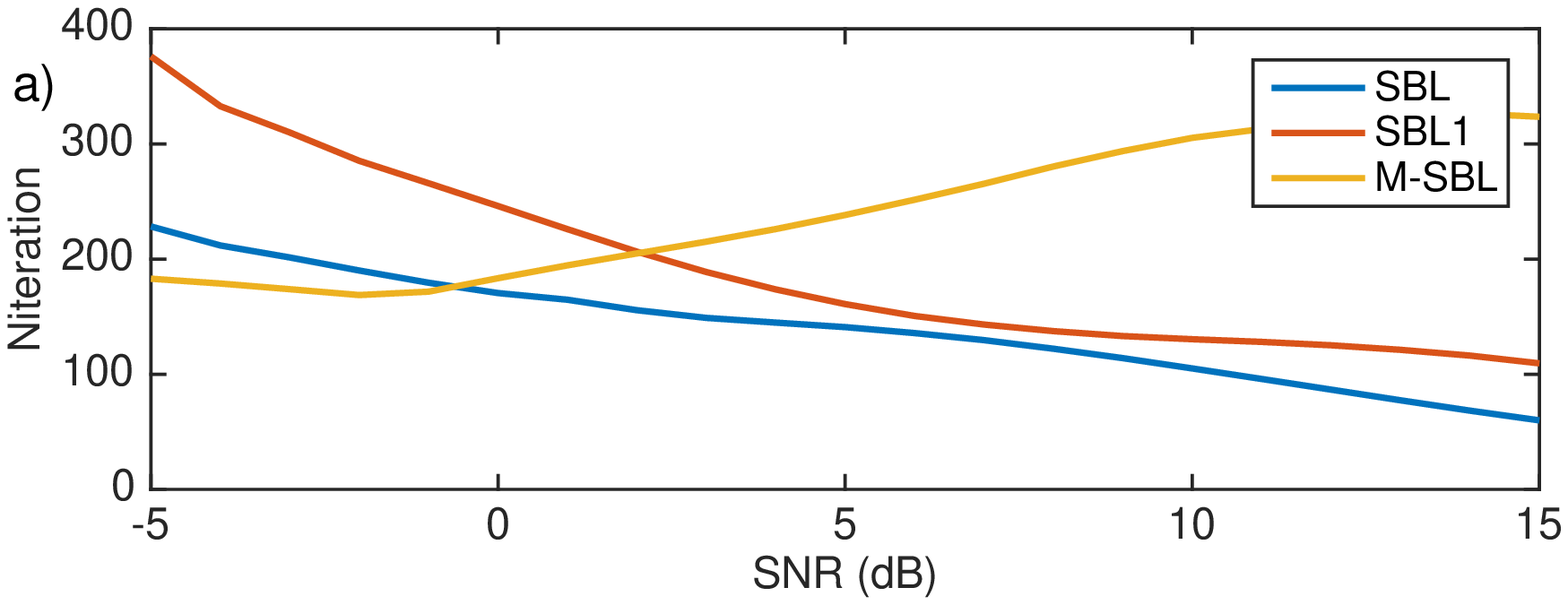}
   \includegraphics[width=0.4\textwidth]{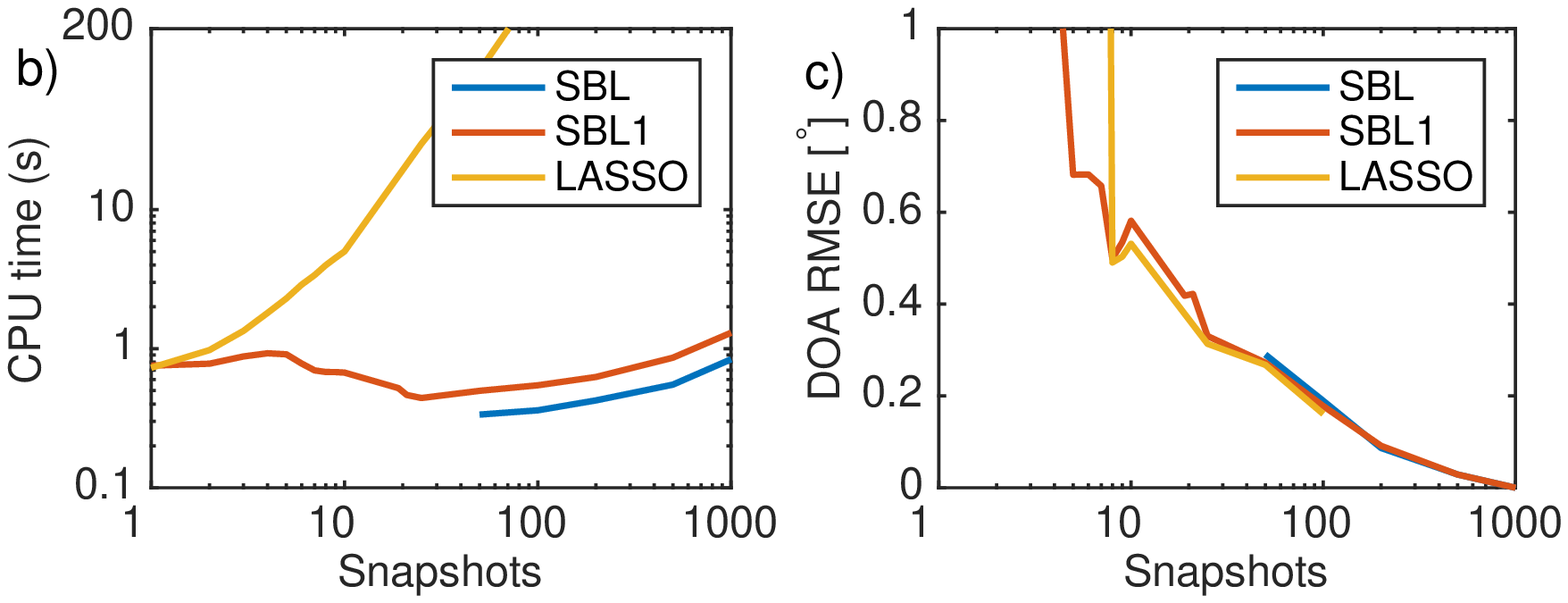}
       \caption{a)  Average number of iterations at each SNR for \ref{eq:gamma-update-EM}, \ref{eq:gamma-update}, and \ref{eq:improved-update} with L=50 snapshots.  
       At array SNR=5 dB, b) average CPU time and c) RMSE for LASSO and \ref{eq:improved-update} vs.\ number of snapshots. 
       All results are an average of 100 Monte Carlo simulations. \label{fig:cpu}}
\end{figure}

The average  number of iterations for \ref{eq:improved-update} and \ref{eq:gamma-update} decreases with SNR but increases for \ref{eq:gamma-update-EM} (Fig.\ \ref{fig:cpu}a). 
%
For \ref{eq:improved-update} and \ref{eq:gamma-update} the CPU time (Macbook Pro 2014) is nearly constant with number of snapshots (Fig.\ \ref{fig:cpu}b). The number of estimated parameters ($\Vecgamma, \sigma^2$)  is independent on the number of snapshots, but increasing the number of snapshots improves the estimation accuracy (lower RMSE). 
%
Contrarily, for LASSO the number of degrees of freedom in ${\Mat{X}}$ increases as  do CPU time with number of snapshots increases.

\section{Conclusions}
A sparse Bayesian learning (SBL) algorithm is derived for high-resolution DOA estimation from multi-snapshot complex-valued array data. The algorithm uses  evidence maximization based on derivatives to estimate the source powers and the noise variance. The method uses the estimated source power at each potential DOA as a proxy for an active DOA promoting sparse reconstruction.

Simulations indicate that the proposed SBL algorithm is a factor of 2 faster than established EM approaches at the same estimation accuracy.
 Increasing the number of snapshots improves the estimation accuracy while the computational effort is nearly independent of snapshots. 
 
\bibliographystyle{unsrt}
\bibliography{CSlasso}

\end{document}